\documentclass[12pt]{amsart}

\usepackage{amssymb}

\newtheorem{proposition}{Proposition}

\newtheorem{corollary}{Corollary}

\newtheorem{fact}{FACT}

\begin{document}

\title[Cyclic cohomology of von Neumann algebras]{Vanishing
of the cyclic cohomology of infinite von Neumann algebras}

\author{Ricardo Bianconi}

\address{Department of Mathematics, IME-USP,
Caixa Postal 66281, CEP 05311-970, S\~ao Paulo, SP,
Brazil}

\email{bianconi@ime.usp.br}

\thanks{\indent
Typeset in AMS-\LaTeX.}

\begin{abstract}
We prove that if $A$ is an infinite von Neumann algebra (i. e., the
identity can be decomposed as a sum of a sequence of pairwise disjoint projections,
all equivalent to the identity)
then the cyclic cohomology of $A$ vanishes. We show that the method of
the proof applies to certain algebras of infinite matrices.
\end{abstract}

\maketitle

\section{Introduction}

Cyclic cohomology is a theory introduced by Alain Connes in \cite{connes84}
with the intent to do the analogue of de Rham cohomology in the context of nocommutative
operator algebras (also independently discovered by B. Tsygan \cite{tsygan83}
and J.-L. Loday and D. Quillen \cite{loday83} in the context of homology
of Lie algebras).

In the paper \cite[Part II, pp 310-360]{connes84}, A. Connes develops the theory from scratch,
and the result we prove here was an effort to make precise the idea
sketched in the end of page 319 and the beginning of the page 320, in which he
indicates the proof of this fact for the algebra of infinite
matrices introduced in \cite{karoubi71}. This we do explicitly in Corollary 2 below.
More precisely, we show by an elementary argument that all the cyclic cohomology
groups of an infinite von Neumann algebra vanish. We recall that an infinite von Neumann
algebra is one which have two (and therefore a parwise disjoint sequence)
disjoint projections whose sum is the identity.

This implies in particular
that if $H$ is an infinite dimensional Hilbert space, then the
cyclic cohomology of the algebra $L(H)$ of bounded operators in $H$ has
trivial cyclic cohomology (as observed by M. Wodzicki in \cite{wodzicki88}, in the
context of homology).
On the other hand, since $H^0_{\lambda}(A)$ is the set of tracial funcionals in
$A$, if the algebra $A$ has a finite trace, then the cohomology does not
vanish.

Since a von Neumann algebra $M$ can be decomposed uniquely as a direct sum
of subalgebras $M_I$, $M_{II_1}$, $M_{II_{\infty}}$ and $M_{III}$, of types
$I$, $II_1$, $II_{\infty}$ and $III$, respectively (see \cite[chapter 2]{sakai71}),
and $M_{II_{\infty}}$, $M_{III}$ and the infinite part of $M_I$ are infinite,
then the cyclic cohomology of $M$ is determined by the summands $M_{II_1}$
(for which there are finite tracial states) and the finite dimensional
parts of $M_I$ (which are Morita equivalent to $\mathbb C$) and so they have nontrivial
cyclic cohomology.

In the following section we prove the main result and derive some corollaries.
We assume that the reader is familiar with \cite[Part II, pp 310-360]{connes84}.

The author wishes to thank Prof. S. T. Melo for helpful discussions. This work
was done while the author recieved research grant Nr. 301238/91-0 from
CNPq (Brazil).

\section{The main result}

The main ingredient we use is the following result by A. Connes.

\begin{fact}\rm
Let $A$ be a unital $\mathbb{C}$-algebra.
Suppose that there is a homomorphism $\rho:A\to A$ and a matrix $V\in{\mathbb{M}}_2(A)$,
such that
\begin{equation*} 
V
\left[
\begin{array}{ll}
a & 0 \\
0 & \rho(a)
\end{array}
\right]
V^{-1}=
\left[
\begin{array}{ll}
0 & 0 \\
0 & \rho(a)
\end{array}
\right].
\end{equation*} 
Then the cyclic cohomology of $A$ vanishes. (See
\cite[Proposition 5.2, page 319]{connes84}.)
\end{fact}

Now, some basic facts about von Neumann algebras.

\begin{fact}\rm
If $A$ is an infinite von Neumann algebra, then
there are partial isometries $w_n,v_n\in A$, $n\in\mathbb N$,
satisfying $p_n=v_nw_n$ is a projection, $w_nv_n=1$, $p_np_m=0$ if $n\neq m$,
and $1=\sum_np_n$. (See \cite[Proposition 2.2.4, page 84]{sakai71}.)

We list some properties of the elements $v_n,w_n\in A$:
\begin{enumerate}
\item
$w_np_n=w_n$,
\item
$p_nv_n=v_n$,
\item
$w_nv_m=0$, if $n\neq m$,
\item
we can assume that $\|v_n\|,\|w_n\|\leq 1$.
\end{enumerate}
\end{fact}

\begin{proposition}\rm
If $A$ is an infinite von Neumann algebra, then
the cyclic cohomology of $A$ vanishes.
\end{proposition}

\proof
Let $\rho:A\to A$, $\rho(a)=\sum_{i\geq 0}v_{2i+1}\,a\,w_{2i+1}$

Let $\theta_1,\theta_2,\theta_3\in A$ be the elements
\begin{equation*} \theta_1=\sum_{i\geq 1}p_i=1-p_0,\end{equation*} 
\begin{equation*} \theta_2=p_0+\sum_{i\geq 1}(v_{2i}w_{2i-1}+v_{2i-1}w_{2i})\end{equation*} 
\begin{equation*} \theta_3=\sum_{i\geq 0}(v_{2i}w_{2i+1}+v_{2i+1}w_{2i})\end{equation*} 

All these sums converge in the strong operator topology and $\|\rho(a)\|\leq \|a\|$.

Observe that if $A$ acts as an algebra of operators in the Hilbert space $H$, then
$\theta_2$ maps $p_0 H$ onto itself and maps $p_{2i} H$ and $p_{2i+1} H$ ($i>0$)
onto each other. Also, $\theta_3$ maps $p_{2i} H$ and $p_{2i+1} H$ ($i\geq 0$)
onto each other.

We have the following relations:
$\theta_1^2=1-p_0$, $\theta_2^2=\theta_3^2=1$, $w_0\theta_1=\theta_1v_0=0$,
 and
\begin{equation*} \theta_1\rho(a)\theta_1=\rho(a)\end{equation*} 
\begin{equation*} \theta_2(v_0aw_0+\rho(a))\theta_2=\sum_{i\geq 0}v_{2i}\,a\,w_{2i}\end{equation*} 
\begin{equation*} \theta_3(\sum_{i\geq 0}v_{2i}\,a\,w_{2i})\theta_3=\rho(a)\end{equation*}

We define the matrices
$X, Y, Z\in {\mathbb{M}}_2(A)$ as follows.

\[
X=\left(
\begin{array}{ll}
0 & w_0 \\
v_0 & \theta_1
\end{array}
\right),\
Y=\left(
\begin{array}{ll}
1 & 0 \\
0 & \theta_2
\end{array}
\right),\
Z=\left(
\begin{array}{ll}
1 & 0 \\
0 & \theta_3
\end{array}
\right).
\] 

Then $X$, $Y$ and $Z$ are invertible and
$X^{-1}=X$, $Y^{-1}=Y$, $Z^{-1}=Z$, and if $V=ZYX$, then $V^{-1}=XYZ$,
and

\[
V
\left[
\begin{array}{ll}
a & 0 \\
0 & \rho(a)
\end{array}
\right]
V^{-1}=
\left[
\begin{array}{ll}
0 & 0 \\
0 & \rho(a)
\end{array}
\right].
\] 

By Fact 1, this implies the vanishing of the cyclic cohomology of $A$.\qed

As corollaries to the proof, we have:

\begin{corollary}\rm
The quotient algebra of an infinite von Neumann algebra by a norm closed
ideal has vanishing cyclic cohomology. In particular,
the Calkin algebra $Q=B(H)/K(H)$ has zero cyclic cohomology,
for $H$ an infinite dimensional
Hilbert space.
\end{corollary}

\proof
If $u,v\in A$, $uv=1$ and $vu=p$, a projection, and $I$ a proper norm closed ideal,
then $p\not\in I$. Therefore the computations done in the proof of the proposition
goes through the quotient.\qed

Let $C$ be the algebra  
introduced in \cite{karoubi71}, of infinite matrices
$(a_{i,j})_{i,j\in\mathbb N}$ such that
\begin{itemize}
\item
the set $\{a_{i,j}:i,j\in\mathbb N\}$ is finite,
\item
the number of nonzero entries in each line and each collumn is bounded.
\end{itemize}

Let $D$ be the algebra  of infinite matrices $(a_{i,j})$ such that
the number of nonzero entries in each line and each collumn is finite.

\begin{corollary}\rm
The cyclic cohomology groups of $C$ and of $D$ vanish.
\end{corollary}

\proof
Repeat the proof of the proposition but now with the sums defined componentwise
(therefore finite in each component). We only need to define the elements
$v_n$ and $w_n$. Let $I_k\subset\mathbb N$, $k\in\mathbb N$, be pairwise disjoint infinite
sets, such that $\mathbb{N}=\bigcup_{k\geq 0}I_k$. Enumerate each $I_k$ as
$I_k=\{i_{k,n}:n\geq 0\}$ (for instance in increasing order). We define the matrices
$v_n=(a^n_{i,j})_{i,j\in\mathbb N}$ (and $w_n$ as the transpose of $v_n$)
as
\begin{equation*} 
a^n_{i,j}=
\left\{
\begin{array}{ll}
1, & \mbox{if $i\in I_n$ and $i=i_{n,j}$,} \\
0, & \mbox{otherwise.}
\end{array}
\right.
\end{equation*} 

Then the entry $i,j$ of $w_nv_n$ is $\sum_{k\geq 0}a^n_{k,i}a^n_{k,j}$,
which is 1 if $i=j$ and 0 otherwise, that is, $w_nv_n=1$. On the other
hand, the entry $i,j$ of $p_n=v_nw_n$ is $\sum_{k\geq 0}a^n_{i,k}a^n_{j,k}$,
which is 1 if $i=j\in I_n$ ($i=j=i_{n,k}\in I_n$, for a unique $k\geq 0$)
and 0 otherwise. Clearly $1=\sum_np_n$, and all these matrices belong to
the algebras $C$ and $D$.
\qed

\bibliographystyle{amsplain}

\end{document}